\documentclass[letter, 10pt, conference]{cssconf}
\IEEEoverridecommandlockouts

\usepackage{alltt,amssymb,latexsym,amsmath,amsthm2,amscd,graphics,epsfig,psfrag,subfigure}

\usepackage{amsfonts,diagrams}
\usepackage{graphicx}
\usepackage{psfrag}

    \theoremstyle{definition}

    \newtheorem{prb}{Problem}[section]

    \newtheorem{lem}{Lemma}[section]

    \newenvironment{rem}{\noindent \textbf{Remark.}}{\hfill{$\bullet$}}

    \newenvironment{pf}{\noindent \textbf{Proof.}}{\hfill{$\blacksquare$}}

\newcommand{\Rbb}{\mathbb{R}}

\newcommand{\Sbb}{\mathbb{S}}

\newcommand{\Lcal}{\mathcal{L}}
\newcommand{\Ical}{\mathcal{I}}

\newcommand{\Ccal}{\mathcal{C}}
\newcommand{\Jcal}{\mathcal{J}}
\newcommand{\Fcal}{\mathcal{F}}

\newcommand{\Ncal}{\mathcal{N}}
\newcommand{\Dcal}{\mathcal{D}}

\newcommand{\Pcal}{\mathcal{P}}
\newcommand{\Qcal}{\mathcal{Q}}
\newcommand{\Rcal}{\mathcal{R}}

\newcommand{\lambdab}{\boldsymbol{\lambda}}

\newcommand{\mub}{\boldsymbol{\mu}}
\newcommand{\etab}{\boldsymbol{\eta}}

\newcommand{\xib}{\boldsymbol{\xi}}
\newcommand{\omegab}{\boldsymbol{\omega}}

\newcommand{\taub}{\boldsymbol{\tau}}

\newcommand{\deltab}{\boldsymbol{\delta}}

\newcommand{\ub}{\mathbf{u}}
\newcommand{\Xb}{\mathbf{X}}
\newcommand{\Yb}{\mathbf{Y}}

\newcommand{\Zb}{\mathbf{Z}}

\newcommand{\qb}{\mathbf{q}}

\newcommand{\Fb}{\mathbf{F}}

\newcommand{\vb}{\mathbf{v}}

\newcommand{\Wb}{\mathbf{W}}

\newcommand{\Ssf}{\textsf{S}}
\newcommand{\Qsf}{\textsf{Q}}
\newcommand{\Tsf}{\textsf{T}}

\newcommand{\Gsf}{\textsf{G}}

\newcommand{\SEsf}{\textsf{SE}}

\newcommand{\drm}{{\rm d}}

\newcommand{\be}{\begin{eqnarray}}
\newcommand{\ee}{\end{eqnarray}}

\newcommand{\nn}{\nonumber}

\begin{document}






\pagestyle{empty}

\title{\LARGE \bf
Optimal Control of Underactuated Nonholonomic Mechanical Systems}

\author{I. Hussein\thanks{Islam Hussein is a Post-Doctoral Fellow at the Coordinated Science Laboratory,  University of Illinois,
Urbana-Champaign, ihussein@uiuc.edu.} \and A. Bloch\thanks{Anthony
Bloch is Professor of Mathematics at the University of Michigan, Ann
Arbor, abloch@umich.edu.}}
\date{}
\maketitle

\begin{abstract}

In this paper we use an affine connection formulation to study an
optimal control problem for a class of nonholonomic, under-actuated
mechanical systems. In particular, we aim at minimizing the
norm-squared of the control input to move the system from an initial
to a terminal state. We consider systems evolving on general
manifolds. The class of nonholonomic systems we study in this paper
includes, in particular, wheeled-type vehicles, which are important
for many robotic locomotion systems. The two special aspects of this
optimal control problem are the nonholonomic constraints and
under-actuation. Nonholonomic constraints restrict the evolution of
the system to a distribution on the manifold. The nonholonomic
connection is used to express the constrained equations of motion.
Furthermore, it is used to take variations of the cost functional.
Many robotic systems are under-actuated since control inputs are
usually applied through the robot's internal configuration space
only. While we do not consider symmetries with respect to group
actions in this paper, the fact that the system is under-actuated is
taken into account in our problem formulation. This allows one to
compute reaction forces due to any inputs applied in directions
orthogonal to the constraint distribution. We illustrate our ideas
by considering a simple example on a three-dimensional manifold.

\end{abstract}

\section{Introduction}\label{eq:Intro}

In this paper we use the theory of affine connections to study force
minimizing optimal control problems for a large class of
nonholonomic under-actuated mechanical systems. Mechanical systems
considered in this paper may be nonlinear and evolving on algebraic
(for holonomically constrained systems) and/or abstract manifolds
such Lie groups (in particular, the group of rigid body motions in
three dimensional space, $\SEsf(3)$, and its subgroups). The class
of nonholonomic systems we study in this paper includes, in
particular, any wheeled-type vehicle, such as robots on wheels and
or tracks. The fact that most of these robotic systems apply torques
and forces internal to the system, which makes these system move in
an undulatory fashion (see \cite{ostrowski:99} and references
therein for more on undulatory locomotion), without the application
of any external forces, makes the system under-actuated. In fact,
control inputs that are applied through the shape space in the
absence of any control authority through the group space (that is,
the \emph{fiber}.) Hence, including under-actuated systems in our
study is crucial in covering a wide range of robotic applications.

Nonholonomic mechanical control systems have a long and complex
history which is described in, for example, \cite{borisov:02},
\cite{bloch:03} (in particular, Chapter 5) and \cite{bullo:04}. Of
much interest in the present work are the recent developments that
utilize a geometric approach \cite{bloch:96,bloch:03} and, in
particular, the theory of affine connections
\cite{Lewis:00,bullo:04}. These methods offer a coordinate-free
differential approach to mechanics and control that avoids many of
the issues that arise in classical mechanics such singularity and
change of coordinates, complexity of notation and the lack of a
geometric picture. For more on differential-geometric mechanics and
its use in the context of dynamics and control, we refer the reader
to \cite{marsden:99,bloch:03,bullo:04}. For the treatment of
under-actuated systems using affine connections, we refer the reader
to \cite{hussein:ACC05UACT}.

Aside from \cite{koon:97,koon:97PhD}, previous results usually treat
kinematic systems that usually aim at minimizing energy. In this
paper the cost function is the square of the norm of the total
applied control. We treat second order (i.e., dynamic) nonholonomic
systems and allow for under-actuation. As will be seen in this
paper, the set of necessary optimality conditions are
coordinate-free and generic for a large class of nonholonomic
mechanical systems. Given problem-specific data, one can specialize
the result to the specific problem at hand. This process can be
automated using symbolic manipulation packages such as
Mathematica$^\circledR$ and toolboxes such as those introduced in
\cite{Bullo:04supp}\footnote{These packages are available online for
which a reference is provided in \cite{Bullo:04supp}.}.

While most of the systems appearing in robotics naturally posses
symmetries with respect to a group action, which leads to the
reduced equations of motion for the system, in this work we provide
a framework for treating nonholonomic systems in the context of
optimal control using the theory of affine connections, regardless
of the presence of any symmetries. In the case where symmetries do
exist, one can usually do more by utilizing the structure of the
equations of motion as done in \cite{koon:97,koon:97PhD}. In
\cite{koon:97,koon:97PhD}, however, the authors use the momentum
equation form of the reduced equations of motion \cite{bloch:96}.
The problem of optimally controlling systems with symmetry will be
treated in a future paper. In particular, we are interested in
understanding how results based on an affine connection approach and
Lagrange's multiplier method relate to results based on the momentum
equation form that appear in \cite{koon:97,koon:97PhD}. For more on
systems with symmetry we refer the reader to, for example,
\cite{ostrowski:99,bloch:96,marsden:99,cendra:01,cendra:01b} and
references therein.

The paper is arranged as follows. In Section \ref{sec:review}, we
briefly describe how nonholonomic mechanical systems are treated
using the theory of affine connections. We also state the
relationship of this approach to the Lagrange-d'Alembert equations
of motion for nonholonomic systems. In Section
\ref{sec:OptimalControl}, we introduce the optimal control problem
and derive the necessary optimality conditions using the theory of
affine connections. In Section \ref{sec:examples}, we use the
vertical coin (equivalently, the inline or ice skate) as a simple
example to illustrate how to perform the computations. Finally, in
Section \ref{sec:conclusion}, we summarize our results and describe
areas of current and future research.

\section{Review of Affine Differential Geometry and Nonholonomic Systems}\label{sec:review}

\subsection{Riemannian Manifolds and Affine Connections}

In this section we give brief definitions of the various objects
from affine connection theory that are essential to this paper. For
more complete studies, we refer the reader to the
mathematically-oriented text \cite{boothby:75} or the more
mechanically-oriented text \cite{bullo:04}.

Let $\Qsf$ be a smooth ($\mathcal{C}^{\infty}$) \emph{Riemannian
manifold} with the \emph{Riemannian metric} defined by
$g_{\qb}:\Tsf_{\qb}\Qsf\times\Tsf_{\qb}\Qsf\rightarrow\Rbb$ at some
point $\qb \in \Qsf$, where $\Tsf\Qsf=\cup_{\qb}\Tsf_{\qb}\Qsf$ is
the \emph{tangent bundle} of all tangent space $\Tsf_{\qb}\Qsf$ at
all points $\qb\in\Qsf$. Thus the length of a tangent vector
$\vb_{\qb} \in \Tsf_{\qb}\Qsf$ is given by
$\sqrt{g_{\qb}(\vb_{\qb},\vb_{\qb})}$.

A \emph{Riemannian connection} on $\Qsf$, denoted $\nabla$, is a
mapping that assigns to any two smooth vector fields $\Xb$ and $\Yb$
on $\Qsf$ a new vector field, $\nabla_{\Xb}\Yb$. For the properties
of $\nabla$, we refer the reader to
\cite{boothby:75,Camarinha:96,bloch:03}. The operator
$\nabla_{\Xb}$, which assigns to every vector field $\Yb$ the vector
field $\nabla_{\Xb}\Yb$, is called the \emph{covariant derivative of
$\Yb$ with respect to $\Xb$}.

The \emph{Lie bracket} of the vector fields $\Xb$ and $\Yb$ will be
denoted by $[\Xb,\Yb]$ and is defined by the identity:
$[\Xb,\Yb]f=\Xb(\Yb f)-\Yb(\Xb f)$, for all $\Ccal^{\infty}$
functions $f:\Qsf\rightarrow \Rbb$. Given vector fields $\Xb$, $\Yb$
and $\Zb$ on $\Qsf$, define the vector field
$\Rcal\left(\Xb,\Yb\right)\Zb$ by the identity
\begin{eqnarray}\label{eq:CurvatureTensorDefinition}
\Rcal\left(\Xb,\Yb\right)\Zb=\nabla_{\Xb}\nabla_{\Yb}\Zb-\nabla_{\Yb}\nabla_{\Xb}\Zb-\nabla_{[\Xb,\Yb]}\Zb.
\end{eqnarray}
$\Rcal$ is trilinear in $\Xb$, $\Yb$ and $\Zb$ and is a tensor of
type $(1,3)$, which is called the \emph{curvature tensor} of $\Qsf$.

Finally, we will employ the musical isomorphism
$\sharp_g:\Tsf^*\Qsf\rightarrow\Tsf\Qsf$ (called the ``sharp") and
its inverse $\flat_g:\Tsf\Qsf\rightarrow\Tsf^*\Qsf$ (the ``flat")
associated with the metric $g$ and defined by the relation
$\Tsf^*\Qsf\ni\Yb^{\flat}(\Xb)=g(\Yb,\Xb)$, for all
$\Xb\in\Tsf\Qsf$. The sharp is induced from the definition of the
flat.

\subsection{Nonholonomic Systems and the Constrained Affine
Connection}

In this section we introduce the affine connection viewpoint of
mechanical control systems. The discussion presented here is based
on the material found in \cite{Lewis:00,bullo:04,cortes:02}. Let
$\Qsf$ be a $\Ccal^{\infty}$ $n$-dimensional manifold with the
tangent and cotangent bundles denoted by $\Tsf\Qsf$ and
$\Tsf^*\Qsf$, respectively. An under-actuated constrained simple
mechanical control system is given by the quadruple
$(\Qsf,g,\Fcal,\Dcal)$, where
$g_{\qb}:\Tsf_{\qb}\Qsf\times\Tsf_{\qb}\Qsf\rightarrow\Rbb$ is the
kinetic energy (Riemannian) metric on $\Qsf$ at $\qb\in\Qsf$. The
collection of covectors
$\Fcal=\left\{\Fb^1,\Fb^2,\ldots,\Fb^p\right\}\in\Tsf^*\Qsf$, $p<n$,
is a set of linearly independent 1-forms on $\Qsf$ that represent
the directions of the forces and torques acting on the system given
by \be\label{eq:tau}\taub=\sum_{i=1}^p\tau_i\Fb^i\in\Tsf^*\Qsf.\ee
Hence, the system is underactuated with underactuation degree $n-p$.
The subspace $\Dcal$ is an $(n-m)$-dimensional nonholonomic
distribution on $\Qsf$, where the $m$ constraints are given by
\be\label{eq:constraints}\omegab^i_{\qb}(\vb_{\qb})=0,
~i=1,\ldots,m,\ee where $\omegab^i\in\Tsf^*\Qsf$ are one-forms on
$\Qsf$. In this paper we only consider systems that evolve in a
potential-free environment. The presence of a potential does not
introduce additional theoretical challenges to our treatment and,
hence, we omit it here for the sake of simplicity. In a future
archival version of this work, we will include potentials to study a
wider range of nonholonomic systems.

For a simple mechanical control system without a potential the
Lagrangian $L:\Tsf_{\qb}\Qsf\rightarrow \Rbb$ is given by
\be\label{eq:Lagrangian} L(\qb,\vb)=\frac{1}{2}g(\vb,\vb).\ee The
Lagrange d'Alembert principle then gives the following equations of
motion \be\label{eq:eoms_Lag_Dal}\frac{\drm}{\drm t}\frac{\partial
L}{\partial\dot{\qb}}-\frac{\partial
L}{\partial\qb}=\sum_{i=1}^k\lambda^j\omegab^j_{\qb}+\sum_{i=1}^m\tau_i\Fb^i,\ee
where $\lambda^j$ are Lagrange multipliers such that
$\lambdab=\sum_{j=1}^k\lambda^j\omegab^j$ represents reaction
forces. The system of equations (\ref{eq:eoms_Lag_Dal}) is
equivalently written using the affine connection as
\be\label{eq:eoms}\nabla_{\vb(t)}\vb(t)&=&\lambdab(t)^{\sharp_{g}}+\ub(t)\nn\\
\dot{\qb}(t)&=&\vb(t)\\
\vb(t)&\in&\Dcal_{\qb(t)}\nn,\ee where $\nabla$ is the Levi-Civita
connection compatible with the metric $g$, $\lambdab(t)$ is a
section of $\Dcal^{\perp}$ (the $g$-orthogonal complement of
$\Dcal$) and \be\label{eq:input}
\ub(t)=\taub^{\sharp_{g}}(t)=\sum_{i=1}^p\tau_i(t)\left(\Fb^i\right)^{\sharp_{g}}=\sum_{i=1}^p\tau_i(t)\Yb_i,\ee
where $\Yb_i=\left(\Fb^i\right)^{\sharp_g}$ are the corresponding
input vector fields.

If we define $\Pcal:\Tsf\Qsf\rightarrow\Dcal\subseteq\Tsf\Qsf$ and
$\Qcal:\Tsf\Qsf\rightarrow\Dcal^{\perp}\subseteq\Tsf\Qsf$ to be the
complementary $g$-orthogonal projectors, then the equations
(\ref{eq:eoms}) are equivalently written
as\be\label{eq:eoms2}\bar{\nabla}_{\vb(t)}\vb(t)&=&\Pcal\left(\ub(t)\right)\nn\\
\dot{\qb}(t)&=&\vb(t),\ee where now we only require that the initial
velocity be $\vb(0)\in\Dcal$ to ensure that the flow remains on the
constrained disctribution. The connection $\bar{\nabla}$ is called
the \emph{nonholonomic affine connection} and is given
by\be\label{eq:NonholonomicAffineConnection}
\bar{\nabla}_{\Xb}\Yb&=&\nabla_{\Xb}\Yb+\left(\nabla_{\Xb}\Qcal\right)\left(\Yb\right)\nn\\&=&\Pcal\left(\nabla_{\Xb}\Yb\right)+\nabla_{\Xb}\left(\Qcal(\Yb)\right),\ee
for all $\Xb,\Yb\in\Tsf\Qsf$. Note that
$\bar{\nabla}_{\Xb}\Yb\in\Dcal$ for all $\Yb\in\Dcal$ and
$\Xb\in\Tsf\Qsf$ \cite{bullo:04,Lewis:00}. The constrained connection
also appears in \cite{Vershik:81}. We now give further
properties of the nonholonomic connection $\bar\nabla$, in
particular, how they operate on functions and one-forms.

\begin{lem}\label{lem:NablaBarF} $\bar{\nabla}_{\Xb}f=\nabla_{\Xb}f$ for all
$f\in\Ccal^{\infty}(\Qsf)$.\end{lem}

\begin{pf} This is obvious since for any affine connection $\tilde{\nabla}$ we have $\tilde\nabla_{\Xb}f=\Lcal_{\Xb}f=\Xb(f)$, the
Lie derivative of $f$ with respect to the vector field $\Xb$. This
is true since the Lie derivative $\Lcal$ is independent of the
choice of $\tilde\nabla$.\end{pf}

\begin{lem}\label{lem:NablaBarLambda} For all $\lambdab\in\Tsf^*\Qsf$ we have
\be\bar{\nabla}_{\Xb}\lambdab=\nabla_{\Xb}\lambdab-
\left(\nabla_{\Xb}\Qcal\right)^*\left(\lambdab\right)
,\nn\ee for all $\Xb\in\Tsf\Qsf$, where $*$ denotes the adjoint of a
map. Note here that $\Qcal$ (and $\Pcal$) is a $(1,1)$ tensor and so
is $\bar\nabla_{\Xb}\Qcal$ and its adjoint
$\bar\nabla_{\Xb}\Qcal$.\end{lem}

\begin{pf} Given our knowledge of how $\bar\nabla$ acts on vector
fields (equation (\ref{eq:NonholonomicAffineConnection})) and Lemma
\ref{lem:NablaBarF}, we have
\be\bar\nabla_{\Xb}\lambdab(\Zb)&=&\bar\nabla_{\Xb}\left(\lambdab(\Zb)\right)-\lambdab\left(\bar\nabla_{\Xb}\Zb\right)\nn\\
&=&\nabla_{\Xb}\left(\lambdab\left(\Zb\right)\right)-\lambdab\left(\nabla_{\Xb}\Zb+\left(\nabla_{\Xb}\Qcal\right)\left(\Zb\right)\right)\nn\\
&=&\left(\nabla_{\Xb}\lambdab\right)\left(\Zb\right)+\lambdab\left(\nabla_{\Xb}\Zb\right)-\lambdab\left(\nabla_{\Xb}\Zb\right)\nn\\&&
-\lambdab\left(\left(\nabla_{\Xb}\Qcal\right)(\Zb)\right)\nn\\
&=&\left(\nabla_{\Xb}\lambdab\right)(\Zb)-\left(\left(\nabla_{\Xb}\Qcal\right)^*\left(\lambdab\right)\right)(\Zb)\nn\ee
for all vector fields $\Xb,\Zb\in\Tsf\Qsf$. For the first equality
we used
$\bar\nabla_{\Xb}\left(\lambdab(\Zb)\right)=\left(\bar\nabla_{\Xb}\lambdab\right)\left(\Zb\right)+\lambdab\left(\bar\nabla_{\Xb}\Zb\right)$.\end{pf}

Finally, recall the definition of the curvature tensor $\Rcal$,
which arises naturally in higher order optimal control problems,
associated with an affine connection $\nabla$ given by equation
(\ref{eq:CurvatureTensorDefinition}). Associated with the
nonholonomic affine connection is the nonholonomic curvature tensor,
denoted $\bar{\Rcal}$, that also satisfies equation
(\ref{eq:CurvatureTensorDefinition}) but with $\bar{\nabla}$
replacing $\nabla$ everywhere. We have the following observation for
the nonholonomic curvature tensor $\bar\Rcal$.

\begin{lem} The
nonholonomic curvature tensor satisfies
\begin{align}\bar{\Rcal}(\Xb,\Yb)\Zb&=\Rcal(\Xb,\Yb)\Zb\nn\\&
+\left(\left(\nabla_{\Xb}\nabla_{\Yb}-\nabla_{\Yb}\nabla_{\Xb}-\nabla_{[\Xb,\Yb]}\right)\Qcal\right)(\Zb)\nn\\
&+(\nabla_{\Xb}\Qcal)\left[\left(\nabla_{\Yb}\Qcal\right)(\Zb)\right]-(\nabla_{\Yb}\Qcal)\left[\left(\nabla_{\Xb}\Qcal\right)(\Zb)\right]\nn\end{align}
for all $\Xb,\Yb,\Zb\in\Tsf\Qsf$.\end{lem}

\begin{pf} The proof comes from the definition in equation (\ref{eq:CurvatureTensorDefinition}) of the nonholonomic curvature
tensor in terms of the nonholonomic connection. One then uses the
definition of the nonholonomic connection in equation
(\ref{eq:NonholonomicAffineConnection}) to substitute $\bar\nabla$
in the equation for $\bar\Rcal$ in
(\ref{eq:CurvatureTensorDefinition}). The rest of the proof is
straightforward algebraic operations.\end{pf}


\section{Optimal Control of a Nonholonomic System}\label{sec:OptimalControl}

In this section we introduce the optimal control problem and derive
the necessary optimality conditions. In this paper, we use
Lagrange's multiplier method for constrained problems in the
calculus of variations. We only investigate normal extremals, which
is a reasonable assumption for simple mechanical control systems
that occur in engineering. We also note that, while the system is
controlled through the shape space only and is, hence, inherently,
under-actuated, \emph{a basic assumption is that the system is
controllable}.

\begin{prb}\label{prb:main_OC_problem_affine} Minimize
\be\label{eq:cost_functional}
\Jcal(\taub)=\int_0^T\frac{1}{2}g^E\left(\ub,\ub\right)\drm t\ee
subject to the dynamics given in equation (\ref{eq:eoms2}) and some
initial and terminal conditions $\qb(0)$ and $\qb(T)$, respectively.
The $(0,2)$ tensor $g^E$ is the standard identity metric on
$\Rbb^n$. The state $\qb(T)$ is assumed to be reachable by the
system from $\qb(0)$.\end{prb}

Note that we want to minimize the norm of $\ub$ as opposed to the
norm of its projection $\Pcal(\ub)$. In other words, generally, the
above formulation does not attempt to minimize the constrained
applied torques. Intuitively, one forecasts that no control forces
and torques should be applied in directions that violate the
constraints since these will be squandered by creating only more
reaction forces (that maintain the constraints) with no net useful
motion, or by violating the constraints altogether. For example, for
the \emph{rolling} vertical coin \cite{bloch:96}, if excessive
torque is applied in the rolling direction, the rolling constraint
may be violated. Moreover, the application of any side forces will
not contribute to the net motion of the system due to the strict
no-side-slip (the ``knife edge") constraint. As will be shown below,
it turns out that the control will be constrained to lie in $\Dcal$
as expected, hence, not allowing for violation of the constraints or
the application of unnecessary control.

In undulatory locomotion, which is of main interest in this work, by
definition, we usually require that the control be applied through
the shape (internal configuration) space only. Hence, we need to
impose the constraint that the generalized control vector field in
the group directions be zero. We do this as follows. Let
$\tilde{\Fb}_i$, $i=1,\ldots,n-p$, form an orthogonal set of
co-vector fields complement to the co-vector fields $\Fb_i$,
$i=1,\ldots,p$. Then, define a Lagrange multiplier one-form
$\xib=\sum_{i=1}^{n-p}\xi_i\tilde{\Fb}_i$ such that
\be\label{eq:control_constraints} \xib(\ub)=0.\ee

We only treat the general case when $\Qsf$ is an arbitrary manifold
and specialize the result to trivial fiber bundles in future
publications. By trivial fiber bundles we mean manifolds of the form
$\Qsf=\Gsf\times\Ssf$, where $\Gsf$ is a Lie group that represents
the fiber, or overall configuration of the system, and $\Ssf$ is the
shape space, or internal configuration, of the system. We begin by
forming the appended cost functional
\be\label{eq:variations1}\Jcal&=&\int_0^T\frac{1}{2}g^E(\ub,\ub)+\xib(\ub)+\mub\left(\dot{\qb}-\vb\right)
\nn\\&&+\etab\left(\bar{\nabla}_{\vb}\vb-\Pcal\left(\ub\right)\right)\drm
t,\ee where $\mub,\etab\in\Dcal^*\Qsf$ are Lagrange multipliers.

\begin{rem} Let $\omegab^j$, $j=m+1,\ldots,n$, be the set of one-forms that span $\Dcal^*$. This set of one forms along with
the one-forms $\omega^i$, $i=1,\ldots,m$, form a basis for $\Tsf^*\Qsf$. In the above, we view $\mub\in\Dcal^*$, the co-tangent constraint distribution spanned by the one forms $\omegab^j$ for , as a one-form on $\Qsf$ such that
$\mub:\Dcal\rightarrow\Rbb$. On the other hand, it is
important to note that the derivatives of the velocity vector field
is general not going to be in the constraint distribution $\Dcal$.
Recall that if $\vb\in\Dcal$, then $\bar\nabla_{\Xb}\vb\in\Dcal$ for all $\Xb\in\Tsf\Qsf$.
Hence, the argument of $\etab$ is always in $\Dcal$ and we view
$\Dcal^*\ni\etab:\Dcal\rightarrow\Rbb$. Moreover, observe that $\Pcal^*\omegab=\omegab$ for all $\omegab\in\Dcal^*$,
where $\Pcal^*:\Dcal^*\subseteq\Tsf^*\Qsf\rightarrow\Tsf^*\Qsf$ is the
adjoint of the map $\Pcal$.\end{rem}

Taking 
variations of equation (\ref{eq:variations1}) we
obtain\begin{align}\deltab\Jcal&=\int_0^Tg^E\left(\nabla_\Wb\ub,\ub+\xib^{\sharp_{g^E}}\right)+\mub\left(\frac{\drm}{\drm
t}\Wb-\nabla_{\Wb}\vb\right)\nn\\&+
\etab\left(\nabla_{\Wb}\bar{\nabla}_{\vb}\vb-\left(\nabla_{\Wb}\Pcal\right)\left(\ub\right)-\Pcal\left(\nabla_{\Wb}
\left(\ub\right)\right)\right)\drm t,\nn\end{align} where
$\nabla_{\Wb}$ is the covariant derivative with respect to the
variation vector field $\Wb\in\Tsf\Qsf$ given
by\be\Wb=\frac{\partial\qb(t,\epsilon)}{\partial\epsilon}\bigg|_{\epsilon=0},\ee
with $\qb(t,\epsilon)$ being the one-parameter variation of the
optimal curve $\qb(t)$. In the above expression, we used the fact
that
$\tilde\nabla_{\Xb}\lambdab\left(\Zb\right)=\left(\tilde\nabla_{\Xb}\lambdab\right)\left(\Zb\right)+\lambdab\left(\tilde\nabla_{\Xb}\Zb\right)$,
for any affine connection $\tilde\nabla$, vector fields $\Xb$ and
$\Zb$ and any co-vector field $\lambdab$ (see page 78 in
\cite{Crampin:86}). When $\lambdab$ is a Lagrange multiplier, say
$\mub$ or $\etab$, the terms
$\left(\nabla_{\Wb}\mub\right)(\dot\qb-\vb)$ and
$\left(\nabla_{\Wb}\etab\right)\left(\bar\nabla_{\vb}\vb-\Pcal\left(\ub\right)\right)$
give us the equations of motion (\ref{eq:eoms2}) later when we set
$\delta\Jcal$ equal to zero. Hence, as usually done in optimal
control theory, these two terms can be omitted without affecting the
rest of the derivation.


We now study the term $\nabla_{\Wb}\bar\nabla_{\vb}\vb$. Using the
definition of the constrained connection in equation
(\ref{eq:NonholonomicAffineConnection}) we get
\be\nabla_{\Wb}\bar\nabla_{\vb}\vb=\nabla_{\Wb}\left[\Pcal\left(\nabla_{\vb}\vb\right)+\nabla_{\vb}\left(\Qcal\left(\vb\right)\right)\right].\nn\ee
However, since the dynamics given by equation (\ref{eq:eoms2})
guarantee that $\vb\in\Dcal$ (assuming $\vb(0)\in\Dcal$), we have
$\Qcal(\vb)=0$. Therefore, we
obtain\be\nabla_{\Wb}\bar\nabla_{\vb}\vb=\left(\nabla_{\Wb}\Pcal\right)\left(\nabla_{\vb}\vb\right)+\Pcal\left(\nabla_{\Wb}\nabla_{\vb}\vb\right).\nn\ee
Recall that the curvature tensor, $\Rcal$, associated with the
unconstrained connection $\nabla$, satisfies the identity given by
equation (\ref{eq:CurvatureTensorDefinition}). Since $\Wb$ is
arbitrary, and hence independent of $\qb$, then $[\Wb,\vb]=0$ such
that \cite{milnor:63}
\be\nabla_{\Wb}\nabla_{\vb}\vb=\nabla_{\vb}\nabla_{\Wb}\vb+\Rcal\left(\Wb,\vb\right)\vb.\nn\ee
Finally, we conclude that
\be\nabla_{\Wb}\bar\nabla_{\vb}\vb&=&\left(\nabla_{\Wb}\Pcal\right)\left(\nabla_{\vb}\vb\right)\nn\\
&&+\Pcal\left(\nabla_{\vb}\nabla_{\Wb}\vb
+\Rcal\left(\Wb,\vb\right)\vb\right).\nn\ee Using this identity, we
obtain\begin{align}\deltab\Jcal&=\int_0^Tg^E\left(\bar{\nabla}_{\Wb}\ub,\ub+\xib^{\sharp_{g^E}}\right)+\mub\left(\frac{\drm}{\drm
t}\Wb-\nabla_{\Wb}\vb\right)\nn\\&+
\etab\big(\Pcal\left(\nabla_{\vb}\nabla_{\Wb}\vb
+\Rcal\left(\Wb,\vb\right)\vb\right)\nn\\&+\left(\nabla_{\Wb}\Pcal\right)\left(\nabla_{\vb}\vb-\ub\right)
-\Pcal\left(\nabla_{\Wb} \left(\ub\right)\right)\big)\drm
t.\nn\end{align}


For the two first terms in the argument of $\mub$ and $\etab$, we
integrate by parts and use the fact that
$\Wb(0)=\Wb(T)=\nabla_{\Wb}\vb(0)=\nabla_{\Wb}\vb(T)=0$ to
obtain\begin{align}\label{eq:integrate_parts}
&\int_0^T\mub\left(\frac{\drm}{\drm t}\Wb\right) \drm t =  -\int_0^T\left(\nabla_{\vb}\mub\right)\left(\Wb\right)\drm t \\
&\int_0^T\etab\left(\Pcal\left(\nabla_{\vb}\nabla_{\Wb}\vb\right)\right)\drm
t
=-\int_0^T\nabla_{\vb}\left(\Pcal^*\left(\etab\right)\right)\left(\nabla_{\Wb}\vb\right)\drm
t\nn\\&\hspace{1.4in}=-\int_0^T\nabla_{\vb}\etab\left(\nabla_{\Wb}\vb\right)\drm
t,\nn\end{align} where we made use of the fact that $\Pcal^*\etab=\etab$ and that $\etab\in\Dcal^*$.

Next, we refer the reader to the properties of a curvature tensor
$\tilde\Rcal$ defined in terms of the connection $\tilde\nabla$ and
a metric $\tilde{g}$ found in \cite{doCarmo:92}, Proposition 2.5 on
page 91. From these properties, one can show that the curvature
satisfies
$\tilde{g}\left(\tilde\Rcal\left(\Wb,\vb\right)\vb,\Xb\right)=\tilde{g}\left(\tilde\Rcal\left(\Xb,\vb\right)\vb,\Wb\right)$,
where $\tilde\Rcal$ is the curvature tensor based on a connection
$\tilde\nabla$ that is compatible with $\tilde{g}$ (this is shown in
\cite{Camarinha:96}, for example). Going through the proof of
Proposition 2.5 in \cite{doCarmo:92}, one can see that all the
derivation can be generalized to an arbitrary metric $\hat{g}$ and
not only to $\tilde{g}$. Hence, we have
$\hat{g}\left(\tilde\Rcal\left(\Wb,\vb\right)\vb,\Xb\right)=\hat{g}\left(\tilde\Rcal\left(\Xb,\vb\right)\vb,\Wb\right)$,
for any positive definite metric $\hat{g}$ on $\Qsf$. In the context
of our problem, this gives
\begin{align}\label{eq:identity_Rcal}\etab\left(\Pcal\left(\Rcal\left(\Wb,\vb\right)\vb\right)\right)&=g^E\left(\left(\Pcal^*\etab\right)^{\sharp_{g^E}},
\Rcal\left(\Wb,\vb\right)\vb\right)\nn\\&=g^E\left(\Wb,\Rcal\left(\etab^{\sharp_{g^E}},\vb\right)\vb\right).\end{align}

Moreover, note
that\be\label{eq:dual_map_of_Pcal}\etab\left(\Pcal\left(\nabla_{\Wb}\ub\right)\right)&=&
\left(\Pcal^*\etab\right)\left(\nabla_{\Wb}\ub\right)\nn\\&=&
g^E\left(\etab^{\sharp_{g^E}},\nabla_{\Wb}\ub\right).\ee
Finally, recall that $\Pcal$ is a $(1,1)$ tensor. Hence
$\nabla_{\Wb}\Pcal$ is also a $(1,1)$ tensor (for this check out any
book on differential geometry, though it is explicitly stated in
\cite{Lewis:00}). So
$\left(\nabla_{\Wb}\Pcal\right)^*:\Tsf^*\Qsf\rightarrow\Tsf^*\Qsf$
is the dual of $\nabla_{\Wb}\Pcal$. With this observation we
have\be\label{eq:dual_map_of_Pcal_nabla}\etab\left(\nabla_{\Wb}\Pcal\left(\nabla_{\vb}{\vb}-\ub\right)\right)&=&
\left(\left(\nabla_{\Wb}\Pcal\right)^*\etab\right)\left(\lambdab^{\sharp_{g^E}}\right)
\\&=&g^E\left(\left(\left(\nabla_{\Wb}\Pcal\right)^*\etab\right)^{\sharp_{g^E}},\lambdab^{\sharp_{g^E}}\right)\nn\ee
where we recall from equation (\ref{eq:eoms}) that $\lambdab$ is the
net reaction generalized force co-vector field.

Using equations
(\ref{eq:integrate_parts})-(\ref{eq:dual_map_of_Pcal_nabla}), we
find
that\be\deltab\Jcal&=&\int_0^Tg^E\left(\nabla_{\Wb}\ub,\ub+\xib^{\sharp_{g^E}}-\etab^{\sharp_{g^E}}\right)\nn\\
&&-\nabla_{\vb}\mub\left(\Wb\right)
+\left[\Rcal\left(\etab^{\sharp_{g^E}},\vb\right)\vb\right]^{\flat_{g^E}}\left(\Wb\right)
\nn\\&&+\left(\left(\nabla_{\Wb}\Pcal\right)^*\etab\right)\left(\lambdab^{\sharp_{g^E}}\right)\nn\\&&
+g^E\left(\nabla_{\Wb}\vb,-\mub^{\sharp_{g^E}}-\left(\nabla_{\vb}\etab\right)^{\sharp_{g^E}}\right)\drm
t.\nn\ee Two important points need to be emphasized. First, the term
$\nabla_{\Wb}\Pcal$ involves variations in the configuration
variables only since $\Pcal$ is an operator that depends only on the
configuration variables. That is, while $\Pcal$ acts on the velocity
vector field $\vb_{\qb}$ at the point $\qb$, $\Pcal$ as a projector
depends solely on the point $\qb$. Thus, this term depends on $\Wb$
only.

The second observation we wish to make is that $\nabla_{\Wb}\ub$ or
$\nabla_{\Wb}\vb$ each can be separated into two terms. For example,
for $\nabla_{\Wb}\ub$, the first variation term involves only
variations in the control components $\tau_i$, $i=1,\ldots,m$, while
the second will involve configuration variations $\Wb$ only coming
through variations of the basis vectors, denoted by $\Yb_i$. This
should be realized in order to completely (and rigorously) separate
variations in configuration, velocity and control variables. For
more on this, see for example the discussion in Section II in
\cite{hussein:04CDCb} and the definitions of the operators $B$ and
$\delta$ therein. It turns out that if we simply ignore this
separation step and treat $\nabla_{\Wb}\ub$ and $\nabla_{\Wb}\vb$ as
terms that involve variations in control and velocity variables
only, respectively, and no variations in the configuration
variables, then we end up with exactly the same result. We emphasize
that the separation step is the correct rigorous mathematical
approach, while ignoring it is not rigorous albeit reduces the
number of steps to obtain the \emph{same correct} result.

From the above discussion, we realize that since $\Wb$,
$\nabla_{\Wb}\ub$ and $\nabla_{\Wb}\vb$ are independent variations
and since for a normal extremal we must have $\deltab\Jcal=0$, we
conclude that along the optimal trajectory we
must have\be\label{eq:nec_conds_affine}\ub+\xib^{\sharp_{g^E}}&=&\etab^{\sharp_{g^E}}\nn\\
\nabla_{\vb}\mub&=&\left[\Rcal\left(\etab^{\sharp_{g^E}},\vb\right)\vb\right]^{\flat_{g^E}}
+\etab\left(\left(\nabla\Pcal\right)\lambdab^{\sharp_{g^E}}\right)\nn\\
\nabla_{\vb}\etab&=&-\mub,\ee where
$\flat_{g^E}$ and $\sharp_{g^E}$ are the musical isomorphisms with
respect to the standard metric $g^E$.

\begin{rem} \textbf{(The Unconstrained, Fully Actuated Problem)} In the unconstrained case
 $\Pcal=\Ical:\Tsf\Qsf\rightarrow\Tsf\Qsf$ is simply the
identity map on $\Tsf\Qsf$ for all $\qb\in\Qsf$. Similarly for
$\Pcal^*$, $\Pcal^*=\Ical^*:\Tsf^*\Qsf\rightarrow\Tsf^*\Qsf$ is the
identity map in $\Tsf^*\Qsf$ for all $\qb\in\Qsf$. $\Qcal$ will have
a null space $\Ncal\left(\Qcal_{\qb}\right)=\Tsf_{\qb}\Qsf$ for all
$\qb\in\Qsf$. In the fully actuated case we have $\xib=0$.
Therefore, in the unconstrained, fully-actuated case the necessary conditions
(\ref{eq:nec_conds_affine}) reduce to
\be\ub&=&\etab^{\sharp_{g^E}}\nn\\
\nabla_{\vb}\mub&=&\left[\Rcal\left(\etab^{\sharp_{g^E}},\vb\right)\vb\right]^{\flat_{g^E}}\nn\\
\nabla_{\vb}\etab&=&-\mub.\nn\ee These are precisely the result
obtained in, say, \cite{hussein:04CDCb}.\end{rem}

\section{Example: The Vertical Coin}\label{sec:examples}

\subsection{Constrained Equations of Motion} 
\begin{figure}[!htb]
\centering
    \psfrag{xy}[][]{\Large{$(q_1,q_2)$}}
    \psfrag{q3}[][]{\LARGE{$q_3$}}
    \psfrag{u1}[][]{\LARGE{$u_1$}}
    \psfrag{u2}[][]{\LARGE{$u_2$}}
    \psfrag{J}[][]{\LARGE{$J$}}
    \psfrag{m}[][]{\LARGE{$m$}}
\resizebox{3.25in}{!}{\includegraphics{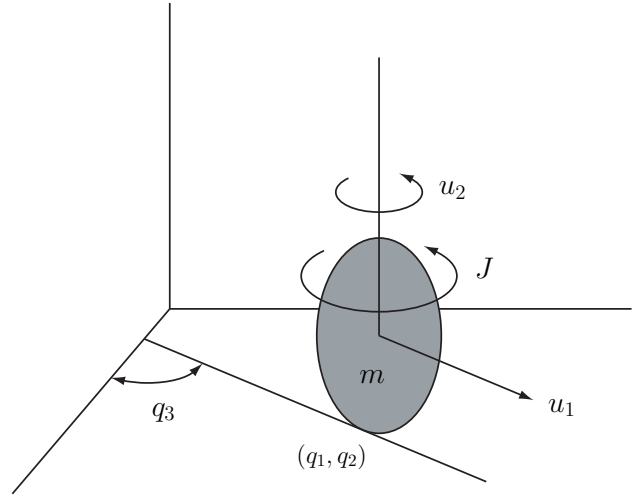}} \caption{The
vertical coin.} \label{Fig:coin}
\end{figure}

We now give an example to illustrate
the above approach and compare the result to traditional methods. In
this section we study optimal control of the vertical coin (i.e., it
can not fall sideways). The system is shown in Figure
\ref{Fig:coin}. The mass of the coin is $m$ and its mass moment of
inertia about the vertical axis is $J$. The position of the point of
contact between the coin and the plane is denoted by $(q_1,q_2)$
while its heading direction is denoted by $q_3$ as shown in the
figure. The configuration $\qb$ is then given by $\qb=(q_1,q_2,q_3)$
and the configuration space is simply $\SEsf(2)$. The control input
is denoted by $u_1$ for the force applied to the center of mass of
the coin and $u_2$ for the torque applied about the vertical axis.

The constraint we have is that the coin can not slip sideways (i.e.,
it satisfies a ``knife-edge" constraint). This constraint is
symbolically expressed in differential form by equation
(\ref{eq:constraints}) with $m=1$ and \be\label{eq:constraints_coin}
\omegab^1_{\qb} = \sin q_3 \drm q_1- \cos q_3 \drm q_2.\ee Hence,
the constraint distribution is given by the span of the vector
fields \be\Xb_1(\qb)&=&\frac{\partial}{\partial
q_3}\nn\\
\Xb_2(\qb)&=&\cos q_3\frac{\partial}{\partial q_1}+\sin
q_3\frac{\partial}{\partial q_2}.\ee The lagrangian for the vertical
coin is given
by\be\label{eq:lagrangian_coin}L(\qb,\dot{\qb})=g(\dot\qb,\dot\qb)=\frac{1}{2}m\left(\dot{q}^2_1+\dot{q}^2_2\right)
+\frac{1}{2}J\dot{q}^2_3\nn\ee and, hence, the components of the
metric $g$ are given by \be\label{eq:metric}g_{11}=g_{22}=m,
~g_{33}=J\ee where all other components are zero.

We denote the unconstrained connection by $\nabla$. Since the metric
is coordinate independent, the unconstrained Christoffel symbols
$\Gamma^i_{jk}$ are all zero. The curvature on $\SEsf(2)$ is
identically zero. (For why curvature is zero on $\SEsf(2)$, see
\cite{hussein:CDC05}.) For the constrained system, one can check
that the Christoffel symbols corresponding to the constrained
connection $\bar{\nabla}$ are given by\be\label{eq:GammaBar}
\bar\Gamma^1_{13}=-\bar\Gamma^2_{23}=\sin (2q_3),
~\bar\Gamma^1_{23}=\bar{\Gamma}^2_{13}=-\cos (2q_3),\ee where all
other Christoffel symbols are zero. 

The external generalized force is given by
\be\label{eq:InputOneForm}\taub=u_1\Fb^1+u_2\Fb^2,\ee where one can
check that the inputs have the directions \be\Fb^1=\cos q_3\drm
q_1+\cos q_3\drm q_2, ~\Fb^2=\drm q_3.\ee The vector fields
$\Yb_i=\left(\Fb^i\right)^{\sharp_g}$ are then given by
\be\Yb_1&=&\frac{\cos q_3}{m}\frac{\partial}{\partial
q_1}+\frac{\sin
q_3}{m}\frac{\partial}{\partial q_2}\nn\\
\Yb_2&=&\frac{1}{J}\frac{\partial}{\partial q_3}\nn\ee such that the
input is given by\be\label{eq:input}\ub=\frac{u_1\cos
q_3}{m}\frac{\partial}{\partial q_1}+\frac{u_1\sin
q_3}{m}\frac{\partial}{\partial
q_2}+\frac{u_2}{J}\frac{\partial}{\partial q_3}.\ee

The projection map $\Pcal:\Tsf\Qsf\rightarrow\Dcal$ is a $(1,1)$
tensor whose elements are computed as follows. Let $\Zb\in\Tsf\Qsf$
be an arbitrary vector field. Then
$\Pcal\left(\Zb\right)=\sum_{i,j}C^{ij}g\left(\Zb,\Xb_i\right)\Xb_j$,
where $C^{ij}$ is the inverse of the matrix
$C_{ij}=g\left(\Xb_i,\Xb_j\right)$, $i,j=1,2$ \cite{cortes:02}. If
we have $\Pcal=\Pcal^i_j\drm q_i \otimes \frac{\partial}{\partial
q_j}$, one can check that the components of $\Pcal$ are given
by\be\label{eq:PBarComponents}\Pcal^1_1&=&\cos^2q_3,
~\Pcal^1_2=\Pcal^2_1=\cos q_3\sin q_3, \nn\\\Pcal^2_2&=&\sin^2q_3,
~\Pcal^3_3=1,\ee where all other components are zero. Hence, we
have\be\Pcal\left(\ub\right)&=&\frac{u_1\cos
q_3}{m}\frac{\partial}{\partial q_1}+\frac{u_1\sin
q_3}{m}\frac{\partial}{\partial
q_2}+\frac{u_2}{J}\frac{\partial}{\partial q_3}=\ub.\nn\ee One could
have anticipated this result since $\ub$ is applied in directions
lying inside the constraint distribution $\Dcal$. This gives the
right hand side of the nonholonomic equation motion given in
equation (\ref{eq:eoms2}).

\begin{rem}\label{rem:UAlongX1} If the force $u_1$ applied to the center of mass was
restricted to be along the $x$-direction, $\Pcal(\ub)$ will not be
equal to $\ub$. In this case, $\ub$ will be projected down to the
constraint distribution by creating a reaction force perpendicular
to the direction of motion to prevent any motion that violates the
knife-edge constraint.\end{rem}

Finally recall that $\bar\nabla_{\vb}\vb$ is given in coordinates by
\be \bar\nabla_{\vb}\vb&=&\sum_{i,j,k=1}^3\left(\ddot
q_i+\bar\Gamma^i_{jk}\dot q_j \dot
q_k\right)\frac{\partial}{\partial q_i}\nn\\
&=&\left[\ddot q_1+\dot q_3\left(\dot q_1\sin (2q_3)-\dot q_2\cos
(2q_3)\right)\right]\frac{\partial}{\partial q_1}\nn\\
&&+\left[\ddot q_2+\dot q_3\left(-\dot q_1\cos (2q_3)-\dot q_2\sin
(2q_3)\right)\right]\frac{\partial}{\partial q_2}\nn\\
&&+\ddot q_3\frac{\partial}{\partial q_3},\nn\ee where we have used
the constrained connection coefficients given in equation
(\ref{eq:GammaBar}). This gives the left hand side of equation
(\ref{eq:eoms2}). Hence, the equations of motion are given by
\be\label{eq:eomsCoin}\ddot{q}_1&=&\left(-\dot q_1\sin(2q_3)+\dot
q_2\cos(2q_3)\right)\dot q_3+\frac{\cos q_3u_1}{m}\nn\\
\ddot{q}_2&=&\left(\dot q_1\cos(2q_3)+\dot q_2\sin(2q_3)\right)\dot
q_3+\frac{\sin q_3u_1}{m}\\
\ddot{q}_3&=&\frac{u_2}{J},\nn\ee which simplify to
\be\label{eq:eomsCoin2}\ddot{q}_1&=&-\dot{q}_2\dot q_3+\frac{\cos q_3u_1}{m}\nn\\
\ddot{q}_2&=&\dot q_1\dot
q_3+\frac{\sin q_3u_1}{m}\\
\ddot{q}_3&=&\frac{u_2}{J}\nn\ee after using the constraints
(\ref{eq:constraints_coin}).

\subsection{Optimality Conditions} We now apply equations
(\ref{eq:nec_conds_affine}) for the vertical coin. First, we need to find a basis for $\Dcal^*$, which we take to be
\be\label{eq:DcalStarBasis} \omegab^2_{\qb}&=&\cos q_3\drm q_1+\sin q_3\drm q_2\nn\\
\omegab^3_{\qb}&=&\drm q_3.\ee Once can check that these are orthogonal to $\omegab^1_{\qb}$. Hence, the lagrange multipliers are given
by \be\label{eq:Mub}\mub&=&\mu^2\omegab^2_{\qb}+\mu^3\omegab^3_{\qb}\nn\\
&=&\mu^2\cos q_3\drm q_1+\mu^2\sin q_3\drm q_2+\mu^3\drm q_3\ee
and \be\label{eq:Etab}\etab&=&\eta^2\omegab^2_{\qb}+\eta^3\omegab^3_{\qb}\nn\\
&=&\eta^2\cos q_3\drm q_1+\eta^2\sin q_3\drm q_2+\eta^3\drm q_3.\ee The under-actuated
direction is given by $\tilde{\Fb}^1=-\sin q_3\drm q_1+\cos q_3\drm
q_2$ such that $\xib=\xi_1\tilde{\Fb}^1$ and
\be\xib^{\sharp_{g^E}}=\frac{-\sin
q_3\xi_1}{m}\frac{\partial}{\partial q_1}+\frac{\cos
q_3\xi_1}{m}\frac{\partial}{\partial q_2}.\nn\ee The first of
equations (\ref{eq:nec_conds_affine}) then gives
\be\label{eq:UNecCond}u_1&=&m\eta^2\nn\\
u_2&=&J\eta^3\\
\xi_1&=&0.\nn\ee In fact, $\xi_1$ is nothing but the generalized reaction force
created by the knife edge constraint in reaction to any forces applied normal to the constraint.

We now obtain the differential equation for $\mub$. The curvature
tensor on $\Qsf$ is identically zero in this case. Moreover, since there are no external forces and torques other than the control inputs, which
are applied in $\Dcal$, then the
constraint reaction forces $\lambdab$ are identically
zero in our example. Generally the reaction forces $\lambdab$ won;t be zero. Finally, $\nabla$ has identically
zero Christofel symbols. Thus, $\nabla_{\vb}$ becomes a simple time
derivative. Hence, the equation for $\mub$ simply gives
\be\dot\mu^2\cos q_3&=&\mu^2\sin q_3\dot q_3\nn\\
\dot\mu^2\sin q_3&=&-\mu^2\cos q_3\dot q_3\nn\\
\dot\mu^3&=&0.\nn \ee Multiplying the first equation by $\cos q_3$ and the second by $\sin q_3$ and summing we obtain
\be\label{eq:MuDotEquation}\dot\mu^2&=&0\nn\\\dot\mu^3&=&0.\ee This gives the $\mub$ differential equation.

Finally, we compute $\nabla_{\vb}\etab=-\mub$. Hence, $\etab$
satisfies the following differential equation
equations
\be\dot\eta^2\cos q_3&=&-\mu^2\cos q_3+\eta^2\dot q_3\sin q_3\nn\\
\dot\eta^2\sin q_3&=&-\mu^2\sin-\eta^2\dot q_3\cos q_3\nn\\
\dot\eta^3&=&-\mu^3.\nn \ee Multiplying the first equation by $\cos q_3$ and the second by $\sin q_3$ and adding both expressions we finally obtain
\be\label{eq:EtaDotEquation}\dot\eta^2&=&-\mu^2\nn\\
\dot\eta^3&=&-\mu^3.\ee

We see that the necessary conditions for this example are particularly simple. This is
due to the fact that the vertical coin has a manifold $\SEsf(2)$ that is isomorphic to $\Rbb^2\times\Sbb^1$, which is a differentially flat
space. This flatness and the positive definiteness and convexity of the cost functional (\ref{eq:cost_functional}) render the problem convex, with simple
linear necessary conditions that can be solved analytically for the global optimal solution. The necessary conditions
are solved to get
\be\label{eq:NecCondSolutions}\mu^2(t)&=&\mu^2_0\nn\\ \mu^3(t)&=&\mu^3_0\nn\\ \eta^2(t)&=&-\mu^2_0t+\eta^2_0\\ \eta^3(t)&=&-\mu^3_0+\eta^3_0,\nn\ee
where $\mu^2_0,\mu^3_0,\eta^2_0,\eta^3_0$ are initial conditions determined from the boundary conditions on the state of the system $\qb$ and $\vb$.
These expressions may then be used to compute the optimal control law from equation (\ref{eq:UNecCond}).

In this paper, we give a simple example for the sake of transparency of the approach. A longer version of this work will include more interesting examples
in the context of optimal robotic locomotion.


\subsection{Verification of Results Using Classical Methods}

In this section we verify the necessary conditions obtained in the
last section, which were the coordinate expressions of equation
(\ref{eq:nec_conds_affine}) for the vertical coin. To do so, we
derive the necessary optimality conditions for the optimal control
problem Problem \ref{prb:main_OC_problem_affine} in coordinates
using equations (\ref{eq:eomsCoin}) for the dynamic constraints. The
cost function in Problem (\ref{prb:main_OC_problem_affine}), in
coordinates,
reads\be\Jcal=\int_0^T\frac{1}{2}\left(\frac{u_1^2}{m^2}+\frac{u_2^2}{J^2}\right)\drm
t.\nn\ee Writing equations (\ref{eq:eomsCoin}) in first order form
\be\label{eq:eomsCoinFirstOrder}\dot{q}_1&=&v_1, ~\dot q_2=v_2,
~\dot q_3=v_3\nn\\\dot v_1&=&-v_2v_3+\frac{\cos q_3u_1}{m}=:f_1(\qb,\vb,\ub)\nn\\
\dot v_2&=&v_1v_3+\frac{\sin q_3u_1}{m}=:f_2(\qb,\vb,\ub)\\
\dot v_3&=&\frac{u_2}{J}=:f_3(\qb,\vb,\ub)\nn\ee the appended cost
function is given by
\be\Jcal&=&\int_0^T\frac{1}{2}\left(\frac{u_1^2}{m^2}+\frac{u_2^2}{J^2}\right)\nn\\&&+\bar\mu^1\left(\dot
q_1-v_1\right)+\bar\mu^2\left(\dot
q_2-v_2\right)+\bar\mu^3\left(\dot q_3-v_3\right)\nn\\
&&+\bar\eta^1\left(\dot v_1-f_1\right)+\bar\eta^2\left(\dot
v_2-f_2\right)+\bar\eta^3\left(\dot v_3-f_3\right) \drm t,\nn\ee
where, following traditional approaches to nonlinear optimal control theory, we take $\bar\mub=\bar\mu^1\drm q_1+\bar\mu^2\drm q_2
+\bar\mu^3\drm q_3$ and $\bar\etab=\bar\eta^1\drm q_1+\bar\eta^2\drm q_2
+\bar\eta^3\drm q_3$.

Taking
variations of the cost functional we
get\be\label{eq:CostVariationsClassical}\delta\Jcal&=&\int_0^T\delta
u_1\left(\frac{u_1}{m^2}-\bar\eta^1\frac{\partial f_1}{\partial
u_1}-\bar\eta^2\frac{\partial f_2}{\partial u_1}\right)\nn\\&&+\delta
u_2\left(\frac{u_2}{J^2}-\bar\eta^3\frac{\partial f_3}{\partial
u_3}\right)\nn\\
&&-\delta q_1 \dot \bar\mu^1-\delta q_2 \dot \bar\mu^2-\delta q_3
\left(\dot{\bar\mu}^3+\bar\eta^1\frac{\partial f_1}{\partial
q_3}+\bar\eta^2\frac{\partial f_2}{\partial q_3}\right)\nn\\
&&+\delta v_1\left(-\bar\mu^1-\dot{\bar\eta}^1-\bar\eta^2v_3\right)\nn\\&&+\delta
v_2\left(-\bar\mu^2-\dot{\bar\eta}^2+\bar\eta^1v_3\right)\nn\\&&+\delta
v_3\left(-\bar\mu^3-\dot{\bar\eta}^3+\bar\eta^1v_2-\bar\eta^2v_1\right).\ee


Setting $\delta\Jcal=0$ in equation
(\ref{eq:CostVariationsClassical}) and by virtue of the independence
of the variations $\delta u_j, ~\delta q_i, ~\delta v_i$, $j=1,2,
~j=1,2,3$, we obtain the following necessary
conditions\be\label{eq:NecCondsClassical}
u_1&=&m\left(\bar\eta^1\cos q_3+\bar\eta^2\sin q_3\right)\nn\\
u_2&=&J\bar\eta^3\nn\\
\dot{\bar\mu}^1&=&0\nn\\
\dot{\bar\mu}^2&=&0\\
\dot{\bar\mu}^3&=&\frac{1}{2}\sin(2q_3)\left[\left(\bar\eta^1\right)^2-\left(\bar\eta^2\right)^2\right]-\bar\eta^1\bar\eta^2\cos(2q_3)\nn\\
\dot{\bar\eta}^1&=&-\bar\mu^1-\bar\eta^2v_3\nn\\
\dot{\bar\eta}^2&=&-\bar\mu^2+\bar\eta^1v_3\nn\\
\dot{\bar\eta}^3&=&-\bar\mu^3+\bar\eta^1v_2-\bar\eta^2v_1.\nn\ee
We need to show that these equations are indeed equivalent to equations (\ref{eq:UNecCond}), (\ref{eq:MuDotEquation}) and (\ref{eq:EtaDotEquation}).
To do that, let $\bar\mub$ and $\bar\etab$ be related to $\mub$ and $\etab$ by \be\bar\mu^1&=&\mu^2\cos q_3, ~\bar\eta^1=\eta^2\cos q_3\nn\\
\bar\mu^2&=&\mu^2\sin q_3, ~\bar\eta^2=\eta^2\sin q_3\\
\bar\mu^3&=&\mu^3, ~\bar\eta^3=\eta^3.\nn\ee
One can easily check that substituting these relationships into equations (\ref{eq:NecCondsClassical}) and after simple algebraic manipulations, we obtain
equations (\ref{eq:UNecCond}), (\ref{eq:MuDotEquation}) and (\ref{eq:EtaDotEquation}).
It is interesting to note the particularly simple form of the necessary conditions obtained using the affine
connection approach when compared to the necessary conditions (\ref{eq:NecCondsClassical}). At first glance, equations (\ref{eq:NecCondsClassical}) may not appear to be solvable, whereas equations (\ref{eq:UNecCond}), (\ref{eq:MuDotEquation}) and (\ref{eq:EtaDotEquation}) are clearly much
easier to study.

\section{Conclusion}\label{sec:conclusion}

In this paper we used the theory of affine connections to study an
optimal control problem for a class of nonholonomic, under-actuated
mechanical systems. The cost function is the norm-squared of the
control input exerted in moving the system from an initial to a
terminal state under the assumption of controllability. We gave a
brief overview of some facts from Riemannian geometry and the use of
the nonholonomic connection to deriving the constrained equations of
motion. We formulated an optimal control problem, where we used the
nonholonomic affine connection together with Lagrange's multiplier
method in the calculus of variations to derive the optimal necessary
conditions. We gave a simple example on a three-dimensional manifold
with a single nonholonomic constraint that captures the main
features of the theoretical result. Future work will focus on the
treatment of nonholonomic systems with symmetry, which naturally
occur in robotic locomotion \cite{ostrowski:99}. In particular, we
are interested in the structure of the resulting optimality
conditions and the possibility of existence of closed form
extremals.

\section*{Acknowledgment}
The authors wish to thank Professor Andrew Lewis for useful
discussions and assistance with the Mathematica packages. We also
wish to acknowledge the scientific input of Dr. Amit Sanyal. The
research of Anthony M. Bloch was supported by NSF grants DMS-030583,
and CMS-0408542.

\bibliographystyle{IEEEtran}
\bibliography{bibliography}


\end{document}